\documentclass[twocolumn]{autart}
\usepackage{indentfirst}
\usepackage{graphicx}
\usepackage{amsmath,amsfonts,amssymb}
\usepackage{natbib}
\setcitestyle{authoryear,round}
\usepackage{nicematrix}
\usepackage{xcolor}
\usepackage{algorithm}
\usepackage{algpseudocode}
\newcommand {\rank} {\mathop{\rm rank}\nolimits}
\allowdisplaybreaks[4]
\numberwithin{equation}{section}
\newtheorem{definition}{Definition}[section]
\newtheorem{theorem}{Theorem}[section]
\newtheorem{lemma}{Lemma}[section]
\newtheorem{problem}{Problem}[section]
\newtheorem{corollary}{Corollary}[section]
\newtheorem{remark}{Remark}[section]
\newtheorem{proposition}{Proposition}[section]
\newtheorem{example}{Example}[section]
\newtheorem{assumption}{Assumption}[section]

\begin{document}

\begin{frontmatter}

\title{Coefficient-level output-feedback stabilization of linear port-Hamiltonian descriptor systems\thanksref{footnoteinfo}}

\thanks[footnoteinfo]{Corresponding author: Juan Zhang.}

\author[XT]{Shuo Shi}\ead{shishuomath@foxmail.com},
\author[XT]{Juan Zhang}\ead{zhangjuan@xtu.edu.cn}

\address[XT]{School of Mathematics and Computational Science, Hunan Key Laboratory for Computation and Simulation in Science and Engineering, Key Laboratory of Intelligent Computing and Information Processing of Ministry of Education,
Xiangtan University, Xiangtan, Hunan, China, 411105}


\begin{abstract}
This paper studies coefficient-level, structure-preserving output-feedback stabilization of linear port-Hamiltonian (pH) descriptor systems.
Existing stabilization conditions generally require explicit pH representations, which may be costly to compute.
We consider descriptor systems for which only the coefficient matrices are available and for which a pH representation is known to exist but is not explicitly given.
For proportional output feedback, we derive coefficient-level conditions that are equivalent to the known solvability criteria in the explicit pH setting.
These conditions ensure that the closed-loop system is regular, impulse-free, asymptotically stable, and remains port-Hamiltonian.
We further extend the framework to proportional-derivative output feedback and enable the assignment of a prescribed dynamical order.
Under the proposed conditions, the proportional gain may be chosen as any symmetric positive definite matrix, and the derivative gain is constructed from coefficient-based decompositions, without computing a pH representation.
\end{abstract}

\end{frontmatter}

\section{Introduction}
Consider the linear time-invariant descriptor system
\begin{align}\label{eqI011}
\left\{ \begin{aligned}
E\dot{x}&=Ax+Bu,\\
y&=Cx,
\end{aligned}\right.
\end{align}
where $x \in \mathbb{R}^n$, $u \in \mathbb{R}^m$, and $y \in \mathbb{R}^m$ denote the state, input, and output vectors, respectively. 
Here, $A, E \in \mathbb{R}^{n\times n}$, $B \in \mathbb{R}^{n\times m}$, and $C \in \mathbb{R}^{m\times n}$ are given coefficient matrices, and the dynamical order of \eqref{eqI011} is $\rank(E)$. 
Descriptor systems combine differential and algebraic components, thereby providing a unified modeling framework for multibody dynamics, chemical and electrical networks, and fluid-mechanical systems
\citep{Dai1989,Duan2010,Eich1998,Feng2017S,Kautsky1989L,Kunkel2006,Mehrmann2023A,Mehrmann2023M,Riaza2008}.

Key properties of descriptor systems include regularity, which is fundamental for well-posedness
\citep{Dai1989,Duan2010,Feng2017S,Kautsky1989L,Xu2006R};
impulse-freeness, which excludes impulsive behavior
\citep{Dai1989,Kautsky1989L};
and asymptotic stability, which is essential in applications
\citep{Duan2010,Xu2006R}.
These properties generally do not hold automatically and must therefore be enforced by suitable feedback control.

Among the available control strategies, static output feedback is particularly appealing because of its structural simplicity and suitability for real-time implementation. 
Extensive research has consequently focused on regularization of descriptor systems by static output feedback
\citep{Gerstner1992S,Gerstner1999L,Gerstner1994I,Chu1998S,Chu1999L,Chu1999I,Nichols2015N,Lewis1990,Shayman1987I}.
Nevertheless, achieving asymptotic stability by output feedback remains highly challenging and is known to be NP-hard
\citep{Blondel1997A,Syrmos1997A},
even for standard state-space systems with $E = I$
\citep{Kailath1980,Zhou1995R}.

\vspace{0.5cm}
Chu and Mehrmann \citep{Chu2024S,Chu2025S} have established solvability conditions for output-feedback stabilization of port-Hamiltonian (pH) descriptor systems when an explicit pH representation is available, while preserving the structure.

The pH descriptor formulation \citep{Beattie2018M} provides an energy-based modeling framework that facilitates modular representation and structure-preserving control of complex physical systems. It has been successfully applied in many areas \citep{Beattie2018M,Mehl2018,Mehrmann2019,Mehrmann2023A,Schaft2023S}.

A remaining challenge is computational: 
when only the coefficient matrices are given and no explicit pH representation is available, one cannot directly apply stabilization criteria that involve pH representations. 
Although a system may admit a pH representation, such representations are generally not unique \citep{Mehrmann2023A}, and computing one may require solving nonlinear matrix constraints or additional realization steps \citep{cherifi2019numerical,van2014port,xu2012direct,Beattie2025A,Chu2025A}.

This difficulty is particularly relevant when the pH structure is only implicit. 
For example, a descriptor system may be known to admit a pH representation if it is positive real and satisfies suitable controllability and observability assumptions \citep{cherifi2024difference,Chu2025A}, yet the specific pH matrices may not be available for controller design. 
In such cases, the question is not to recover a particular representation, but to verify and enforce closed-loop properties directly from the coefficient matrices.

This paper does not aim to determine whether a given system is port-Hamiltonian, nor to compute a pH representation. 
Under the assumption that at least one pH representation exists, we derive solvability conditions and construct feedback gains using only the coefficient matrices, without forming the pH representation matrice.

Although regularization results without explicit pH representations exist \citep{Chu2025R}, structure-preserving output-feedback stabilization has so far been solved mainly for systems with explicit representations \citep{Chu2024S,Chu2025S}. 
This paper removes that limitation by developing coefficient-level tests and constructive designs that avoid constructing any pH representations. 
Our main contributions are:

1. 
We show that solvability of the stabilization problems implies equivalence with the pH representation form studied in \citet{Chu2024S}, thereby clarifying the connection with that framework.

2. 
For proportional output feedback, we reformulate the solvability conditions for explicit pH systems using only the coefficient matrices. 
The proportional gain may be any symmetric positive definite matrix.

3. 
For proportional-derivative output feedback, we provide coefficient-based decompositions and a constructive procedure for selecting the derivative gain so as to assign a prescribed dynamical order, using only the coefficient matrices.

The paper is organized as follows. 
Section \ref{Se2} introduces the problem formulation, key definitions, and preliminary results. 
Section \ref{Se3} presents solvability conditions and gain selection for proportional output feedback without an explicit pH representation.
Section \ref{Se4} extends these results to proportional-derivative output feedback with a prescribed dynamical order. 
The paper concludes with a summary of the main findings.

\textit{\bf Notation.}
Throughout the paper, $I_n$ denotes the $n \times n$ identity matrix. 
For a complex scalar $s \in \mathbb{C}$, its real part is denoted by $\Re(s)$. 
For a real matrix $A$, $\det(A)$, $\rank(A)$, and $A^T$ denote its determinant, rank, and transpose, respectively. 
The notation $A \ge 0$ ($A > 0$) means that $A$ is symmetric positive semidefinite (positive definite). 
A matrix whose orthonormal columns span the right (left) nullspace of $M$ is denoted by $S_{\infty}(M)$ ($T_{\infty}(M)$). 
Finally, the degree of a polynomial $f$ is denoted by $\deg(f)$.

\section{Problem formulation and preliminary results}\label{Se2}
This section introduces the basic notions for descriptor systems, recalls pH descriptor systems, formulates the feedback problems, and states auxiliary lemmas.

We begin by recalling the concepts of regularity, impulse-freeness, and asymptotic stability
\citep{Kunkel2006,Xu2006R}.

\begin{definition}\rm \label{defDuan2010}
The descriptor system \eqref{eqI011}, or equivalently the matrix pair $(E,A)$, is called regular if there exists $s\in \mathbb{C}$ such that $\det(sE-A)\neq0$.
\end{definition}

\begin{definition}\rm
The descriptor system \eqref{eqI011}, or equivalently the matrix pair $(E,A)$, is called regular and impulse-free if $\deg(\det(sE-A))=\rank(E)$.
\end{definition}

\begin{definition}\rm\label{def202412031413}
The descriptor system \eqref{eqI011}, or equivalently the matrix pair $(E,A)$, is called regular and asymptotically stable if all the roots of $\det(sE-A)=0$ have strict negative real parts.
\end{definition}

We next recall the definition of a port-Hamiltonian (pH) descriptor system \citep{Mehrmann2023A}.

\begin{definition}\rm\label{def202508011451}
The linear time-invariant descriptor system \eqref{eqI011} is said to admit a pH representation if its coefficients can be written as
\begin{equation}\label{eq20260329192801}
\begin{bmatrix}
A&B\\
C&0\\
\end{bmatrix}=\begin{bmatrix}
(J-R)Q&G-P\\
(G+P)^TQ&0\\
\end{bmatrix},
\end{equation}
where the given matrix $E$ and matrices $J$, $R$, $Q\in \mathbb{R}^{n\times n}$ and $G$, $P\in \mathbb{R}^{n\times m}$ satisfy
\begin{equation}\label{eq20260329191301}
\begin{aligned}
Q^TE=E^TQ\geq0,~&Q^TJQ=-(Q^TJQ)^T,\\
Q^TRQ\geq0,~&Q^TP=0.
\end{aligned}
\end{equation}
The system \eqref{eqI011} is called a pH descriptor system if it admits at least one such representation.
\end{definition}

Throughout the paper, we work under the following assumption:
    
\begin{assumption}\label{assu:pH}
\rm
The descriptor system \eqref{eqI011} is assumed to admit a pH representation that is not explicitly given; 
that is, there exist matrices $J$, $R$, $Q$, $G$, $P$ satisfying Definition \ref{def202508011451}, yet only the coefficient matrices $E$, $A$, $B$, and $C$ are available.
\end{assumption}

For any chosen pH representation, the quadratic function $\mathcal H(x)=\frac{1}{2}x^TE^TQx$ is called the Hamiltonian and can be interpreted as the stored energy associated with that representation. 
It satisfies the power balance equation \citep{Beattie2018M,Mehrmann2023A},
\begin{equation*}
\frac{d}{dt}\mathcal H(x) = -\begin{bmatrix}
x\\
u\\
\end{bmatrix}^T
\begin{bmatrix}
Q^TRQ&Q^TP\\
P^TQ&0\\
\end{bmatrix}
\begin{bmatrix}
x\\
u\\
\end{bmatrix} + y^Tu
\end{equation*}
along any solution $x$ and for any input $u$. Here, the first term on the right-hand side represents the dissipated energy, whereas the second term represents the supplied energy.

It is well known that pH descriptor systems satisfy the dissipation inequality
\begin{equation*}
\frac{d}{dt}\mathcal H(x) \leq y^Tu,
\end{equation*}
which corresponds to passivity in control theory. 
For this to hold in the case of systems without a feedthrough term, it is necessary that
\begin{equation*}
\begin{bmatrix}
Q^TRQ&Q^TP\\
P^TQ&0\\
\end{bmatrix}\geq 0,
\end{equation*}
which implies $Q^TP=0$ and $Q^TRQ\geq 0$. 
In particular, the uncontrolled system (with $u=0$) is semidissipative.

We now formulate the two output-feedback problems studied in this paper.

\begin{problem}\label{pr20250401143601}
\rm
Consider the descriptor system \eqref{eqI011} under Assumption \ref{assu:pH}.  
Applying the proportional output feedback law $u=-Ky+v$ to the system \eqref{eqI011} yields the closed-loop system
\begin{align}\label{eq2024111601}
\left\{\begin{aligned}
E\dot{x}&=(A-BKC)x+Bv,\\
y&=Cx.
\end{aligned}\right.
\end{align}
The objective is to determine a matrix $K\in\mathbb{R}^{m\times m}$ such that the matrix pair $(E, A - BKC)$ is regular, impulse-free, and asymptotically stable, while the closed-loop system \eqref{eq2024111601} still admits a pH representation, i.e., satisfies Definition \ref{def202508011451}.
\end{problem}

\begin{problem}\label{pr2025072143601}
\rm
Consider the descriptor system \eqref{eqI011} under Assumption \ref{assu:pH}.  
Applying the proportional-derivative output feedback law $u=-Ky-F\dot{y}+v$ to the system \eqref{eqI011} yields the closed-loop system
\begin{align}\label{eq202411160142}
\left\{\begin{aligned}
(E+BFC)\dot{x}&=(A-BKC)x+Bv,\\
y&=Cx.
\end{aligned}\right.
\end{align}
The objective is to find matrices $K$ and $F\in\mathbb{R}^{m\times m}$ such that the pair $(E + BFC, A - BKC)$ is regular, impulse-free, asymptotically stable, has the prescribed dynamical order $\rank(E + BFC) = r$, and admits a pH representation, i.e., satisfies Definition \ref{def202508011451}.
\end{problem}

\begin{remark}
\rm
Whenever we say that a closed-loop system remains port-Hamiltonian, this means only that the corresponding closed-loop coefficient matrices admit at least one pH representation. 
It does not presuppose that a particular open-loop representation has been fixed in advance or explicitly computed.
\end{remark}

The next proposition shows that solvability of either feedback problem forces the matrix $Q$ in every pH representation to be nonsingular.
This observation explains why the underlying solvability question reduces to the explicit pH descriptor framework, while the main remaining issue is representation-free verification and feedback construction.

\begin{proposition}\rm\label{re202510181952}
Consider the descriptor system \eqref{eqI011} under Assumption \ref{assu:pH}.  
If either Problem \ref{pr20250401143601} or Problem \ref{pr2025072143601} is solvable, then

(i)~every pH representation of \eqref{eqI011} has nonsingular $Q$;

(ii)~the relation $Q^TP=0$ implies $P=0$ for every pH representation;

(iii)~for every pH representation, the system \eqref{eqI011} is equivalent to the following pH descriptor system:
\begin{align}\label{eq:Eq-pH}
\left\{\begin{aligned}
(Q^TE)\dot{x}&=(Q^TJQ-Q^TRQ)x+(Q^TG)u,\\
y&=(Q^TG)^Tx,
\end{aligned}\right.
\end{align}
where $Q^TE\geq0$, $(Q^TJQ)^T=-(Q^TJQ)$, $Q^TRQ\geq0$.
\end{proposition}
\textbf{Proof.~~}It suffices to prove the claim for Problem \ref{pr20250401143601}; the other case is analogous.

Choose any pH representation of \eqref{eqI011}. If Problem \ref{pr20250401143601} is solvable, then there exists $K$ such that the closed-loop system \eqref{eq2024111601} is port-Hamiltonian, regular, and asymptotically stable. 
With respect to this representation, the closed-loop system matrix is $(J-(R+(G-P)K(G+P)^T))Q$. Hence,
$$\rank(sE-(J-(R+(G-P)K(G+P)^T))Q)=n$$
for all $s\in \mathbb{C}$ with $\Re(s)\geq0$, by Definition \ref{def202412031413}.
Taking $s=0$ gives
$$\rank((J-(R+(G-P)K(G+P)^T))Q)=n\leq \rank(Q)\leq n.$$
Hence $\rank(Q)=n$, so $Q$ is nonsingular. 
Since the pH representation was chosen arbitrarily, the same conclusion holds for every matrix $Q$. 

Because $Q$ is nonsingular and $Q^TP=0$, it follows that $P=0$. 
Substituting $P=0$ into \eqref{eq20260329192801}, we obtain the reduced representation in part (iii). 
Therefore, both (ii) and (iii) hold for every pH representation.~\hfill$~~\square$

\begin{remark}
\rm
In a general pH descriptor representation, both $E$ and $Q$ may be singular. 
Although Definitio \ref{def202508011451} includes $P$, Proposition \ref{re202510181952} shows that solvability forces $Q$ to be nonsingular and $P=0$. 
Hence, any stabilizable pH descriptor system is equivalent to \eqref{eq:Eq-pH}. 
The stabilization of this form has been studied in \citep{Chu2024S}, which also shows that systems with explicit pH representations can be transformed into it.

The difficulty addressed here is practical: 
under Assumption \ref{assu:pH}, an explicit representation is not available, and applying the results of \citep{Chu2024S} would require computing one, which may be nontrivial. 
This paper therefore does not introduce a new stabilizability class. 
It provides coefficient-level conditions and feedback constructions using only $E$, $A$, $B$, and $C$. 
The proofs are carried out at the coefficient level, so no pH representation is computed. 
The contribution is a representation-free reformulation of the synthesis problem.
\end{remark}

We now state several lemmas that will be used in the analysis.

\begin{lemma}\rm\citep{Dai1989}\label{thDai1989}
The descriptor system \eqref{eqI011} is regular and impulse-free if and only if
$$\rank\begin{bmatrix}
E&0\\
A&E\\
\end{bmatrix}=n+\rank(E).$$
\end{lemma}

\begin{lemma}\rm\citep{Chu2024S}\label{leChu2024S}
Let $E$, $J$, $R\in \mathbb{R}^{n\times n}$ with $E\geq0$, $J=-J^T$, $R\geq0$ and
\vspace{-0.3cm}
$$\begin{bNiceArray}{c|c|c}[nullify-dots]
E& J & R
\end{bNiceArray}=\begin{bNiceArray}{cc|cc|cc}[first-row,first-col,nullify-dots]
&n_1&n_2&n_1&n_2&n_1&n_2\\
n_1&E_{11}&E_{12}&J_{11}&J_{12}&R_{11}&0\\
n_2&E_{12}^T&E_{22}&-J_{12}^T&J_{22}&0&0
\end{bNiceArray},$$
where $R_{11}>0$. Then

(i) The matrix $J-R$ is nonsingular if and only if
\begin{align*}
\rank\begin{bmatrix}
J_{12}\\
J_{22}\\
\end{bmatrix}=n_2.
\end{align*}

(ii) The matrix pair $(E,J-R)$ is regular and asymptotically stable if and only if, for all $s\in \mathbb{C}$ with $\Re(s)=0$,
\begin{align*}
\rank\begin{bmatrix}
J_{12}-sE_{12}\\
J_{22}-sE_{22}\\
\end{bmatrix}=n_2.
\end{align*}
\end{lemma}

\begin{lemma}\rm\label{th2.4}\citep{Chu1998S}
Let $E \in \mathbb{R}^{n\times n}$, $B \in \mathbb{R}^{n\times m}$, and $C \in \mathbb{R}^{p\times n}$. Then there exists $G_0\in \mathbb{R}^{m\times p}$ such that $\rank(E+BG_0C)=r$ if and only if
\begin{align*}
\rank\begin{bmatrix}
E&B\\
\end{bmatrix}+&\rank\begin{bmatrix}
E\\
C\\
\end{bmatrix}-\rank\begin{bmatrix}
E&B\\
C&0\\
\end{bmatrix}\leq r\\
&\leq \min\left\{\rank\begin{bmatrix}
E&B\\
\end{bmatrix},\rank\begin{bmatrix}
E\\
C\\
\end{bmatrix}\right\}.
\end{align*}
\end{lemma}

\section{Reformulation of proportional-feedback stabilization}\label{Se3}

Although the proof invokes an auxiliary pH representation, the conditions and gain selection below are expressed solely in terms of $E$, $A$, $B$, and $C$.

\begin{theorem}\rm \label{eq202507222239}
Consider the descriptor system \eqref{eqI011} under Assumption \ref{assu:pH}.  
Then the following statements are equivalent.

\noindent
(a) There exists $K\in\mathbb{R}^{m\times m}$ for which $(E,A-BKC)$ is regular, impulse-free, and asymptotically stable, and the closed-loop system \eqref{eq2024111601} still admits a pH representation.

\noindent
(b) The following coefficient-level conditions hold:

(i)
\begin{align}\label{eq2024120121250101}
\rank\begin{bmatrix}
E&A\\
0&E\\
0&C\\
\end{bmatrix}=n+\rank(E).
\end{align}

(ii) For all $s\in\mathbb{C}$ with $\Re(s)=0$,
\begin{align}\label{eq20241203143301}
\rank\begin{bmatrix}sE-A&B\end{bmatrix}
=
\rank\begin{bmatrix}sE-A\\
C\end{bmatrix}
=n.
\end{align}
Moreover, whenever (b) holds, $K$ may be chosen as any symmetric positive definite matrix.
\end{theorem}

\textbf{Proof.~~}``(a) $\Longrightarrow$ (b)''
Let $K$ be such that the matrix pair $(E,A-BKC)$ is regular, impulse-free, and asymptotically stable. Then $\rank(sE-(A-BKC))=n$ for all $s$ with $\Re(s)=0$, by Definition \ref{def202412031413}. Hence condition (ii) holds.

Moreover, by Lemma \ref{thDai1989},
\begin{align*}
\rank\begin{bmatrix}
E&0\\
A-BKC&E\\
\end{bmatrix}=n+\rank(E).
\end{align*}
Since
$$
\rank\begin{bmatrix}
E&0\\
A-BKC&E\\
\end{bmatrix}\le n+\rank(E),
$$
the matrix $K$ may be viewed as a variable chosen to maximize this rank. By Lemma \ref{th2.4}, this is equivalent to
\begin{align*}
\min\left\{\rank\begin{bmatrix}
E&0&0\\
A&E&B\\
\end{bmatrix},~\rank\begin{bmatrix}
E&A\\
0&E\\
0&C\\
\end{bmatrix}\right\}=n+\rank(E),
\end{align*}
which shows that condition (i) holds.

``(b) $\Longrightarrow$ (a)''
Assume that conditions (i) and (ii) hold. 
We prove sufficiency by transforming the system into an auxiliary pH representation.

\emph{Step 1}: Reduction to $P=0$. 

Choose a pH representation of \eqref{eqI011}. 
We first show that it induces an equivalent auxiliary pH system with $P=0$.

From condition (ii),
\begin{align*}
\rank\begin{bmatrix}sE-A\\
C\end{bmatrix}
=
\rank\begin{bmatrix}sE-(J-R)Q\\
(G+P)^TQ\end{bmatrix}=n
\end{align*}
for all $s\in \mathbb{C}$ with $\Re(s)=0$.
Taking $s=0$ yields
\begin{align*}
\rank\begin{bmatrix}-(J-R)Q\\
(G+P)^TQ\end{bmatrix}=n\leq \rank(Q)\leq n,
\end{align*}
and hence $\rank(Q)=n$.

Since $Q$ is nonsingular and $Q^TP=0$, we obtain $P=0$. System \eqref{eqI011} is equivalent to the following pH descriptor system:
\begin{align}\label{eq202604052056}
\left\{\begin{aligned}
\widetilde{E}\dot{x}&=(\widetilde{J}-\widetilde{R})x+\widetilde{B}u,\\
y&=\widetilde{B}^Tx,
\end{aligned}\right.
\end{align}
where $\widetilde{E}=Q^TE$, $\widetilde{J}=Q^TJQ$, $\widetilde{R}=Q^TRQ$, $\widetilde{B}=Q^TG$, and $\widetilde{E}\geq0$, $\widetilde{J}^T=-\widetilde{J}$, $\widetilde{R}\geq0$.

\emph{Step 2}: Coefficient conditions in the auxiliary coordinates.

Since $Q$ is nonsingular, conditions (i) and (ii) are equivalent to the following two conditions:

(1)
$$
\rank\begin{bmatrix}
\widetilde{E}&\widetilde{J}-\widetilde{R}\\
0&\widetilde{E}\\
0&\widetilde{B}^T\\
\end{bmatrix}=n+\rank(\widetilde{E}).
$$

(2) For all $s\in \mathbb{C}$ with $\Re(s)=0$,
\begin{align*}
\rank\begin{bmatrix}s\widetilde{E}-(\widetilde{J}-\widetilde{R})&\widetilde{B}\end{bmatrix}
=
\rank\begin{bmatrix}s\widetilde{E}-(\widetilde{J}-\widetilde{R})\\
\widetilde{B}^T\end{bmatrix}=n.
\end{align*}

Therefore, it remains to show that, under conditions (1) and (2), the matrix pair $(\widetilde{E},\widetilde{J}-(\widetilde{R}+\widetilde{B}K\widetilde{B}^T))$ is regular, impulse-free, and asymptotically stable, and that the resulting closed-loop system
\begin{align}\label{eq202507221443}
\left\{\begin{aligned}
\widetilde{E}\dot{x}&=(\widetilde{J}-(\widetilde{R}+\widetilde{B}K\widetilde{B}^T))x+\widetilde{B}v,\\
y&=\widetilde{B}^Tx,
\end{aligned}\right.
\end{align}
also preserves the port-Hamiltonian structure.

\emph{Step 3}: Admissibility of arbitrary positive definite proportional gains.

The system \eqref{eq202604052056} serves only as an auxiliary representation used to prove that conditions (1) and (2) imply the desired closed-loop properties for every $K>0$.

Let $K$ be any symmetric positive definite matrix. 
Using the singular value decomposition, there exist orthogonal matrices $U$ and $V$ such that
\begin{align}\label{eq20241203144401}
\hspace{-0.5cm}
\begin{bNiceArray}{c|c}[nullify-dots]
U\widetilde{R}U^T&U\widetilde{B}V
\end{bNiceArray}=
\begin{bNiceArray}{ccc|cc}[first-row,first-col,nullify-dots]
&n_1&n_2& n_3&n_1&m_1\\
n_1&R_{11}&R_{12}&0&\Sigma_G&0\\
n_2&R_{12}^T&\Sigma_R&0&0&0\\
n_3&0&0&0&0&0
\end{bNiceArray},
\end{align}
where $m_1=m-n_1$, $\Sigma_R>0$, and $\Sigma_G>0$; and let
\begin{equation*}
\begin{aligned}
&\begin{bNiceArray}{c|c}[nullify-dots]
U\widetilde{E}U^T&U\widetilde{J}U^T
\end{bNiceArray}
=
\begin{bNiceArray}{ccc|ccc}[first-row,first-col,nullify-dots]
&n_1&n_2& n_3&n_1&n_2& n_3\\
n_1&E_{11}&E_{12}&E_{13}&J_{11}&J_{12}&J_{13}\\
n_2&E_{12}^T&E_{22}&E_{23}&-J_{12}^T&J_{22}&J_{23}\\
n_3&E_{13}^T&E_{23}^T&E_{33}&-J_{13}^T&J_{23}^T&J_{33}
\end{bNiceArray}.
\end{aligned}
\end{equation*}

Let
$$V^TKV =\begin{bNiceArray}{cc}[first-row,first-col,nullify-dots]
&n_1&m-n_1\\
n_1&K_{11}&K_{12}\\
m-n_1&K_{12}^T&K_{22}
\end{bNiceArray}>0.$$
Then
\begin{align*}
&U\widetilde{R}U^T+U\widetilde{B}VV^TKVV^T\widetilde{B}^TU^T\\
=&\begin{bmatrix}
R_{11}+\Sigma_GK_{11}\Sigma_G&R_{12}&0\\
R_{12}^T&\Sigma_R&0\\
0&0&0\\
\end{bmatrix}\geq 0.
\end{align*}
Since $K>0$ implies $K_{11}>0$, the Schur complement \citep{zhan2013matrix} gives
\begin{equation}\label{eq20241203144701}
\begin{bmatrix}
R_{11}+\Sigma_GK_{11}\Sigma_G&R_{12}\\
R_{12}^T&\Sigma_R
\end{bmatrix}>0.
\end{equation}
Moreover,
\begin{flalign}\label{eq20241203150201}
\rank(\widetilde{R}+\widetilde{B}K\widetilde{B}^T)=\rank\begin{bmatrix}\widetilde{R}&\widetilde{B}\end{bmatrix}.
\end{flalign}

From condition (2), we have
\begin{align}\label{eq20241203160001}
\rank\begin{bmatrix}
sE_{13}-J_{13}\\
sE_{23}-J_{23}\\
sE_{33}-J_{33}\\
\end{bmatrix}=n_3.
\end{align}
Hence, by \eqref{eq20241203144701}, \eqref{eq20241203160001}, and Lemma \ref{leChu2024S}, the matrix pair $(\widetilde{E},\widetilde{J}-(\widetilde{R}+\widetilde{B}K\widetilde{B}^T))$ is regular and asymptotically stable, and the resulting closed-loop system \eqref{eq202507221443} still admits a pH representation.

Finally, there exists an orthogonal matrix $\widetilde{U}$ such that
\begin{equation}\label{eq20241203160801}
\begin{aligned}
&\begin{bNiceArray}{c|c}[nullify-dots]
\widetilde{U}\widetilde{E}\widetilde{U}^T & \widetilde{U}(\widetilde{R}+\widetilde{B}K\widetilde{B}^T)\widetilde{U}^T\\
\end{bNiceArray}\\
=&\begin{bNiceArray}{ccc|ccc}[first-row,first-col,nullify-dots]
&\widetilde{n}_1&\widetilde{n}_2& \widetilde{n}_3&\widetilde{n}_1&\widetilde{n}_2& \widetilde{n}_3\\
\widetilde{n}_1&\Sigma_E&0&0&\widetilde{R}_{11}&\widetilde{R}_{12}&0\\
\widetilde{n}_2&0&0&0&\widetilde{R}_{12}^T&\widetilde{\Sigma}_R&0\\
\widetilde{n}_3&0&0&0&0&0&0\\
\end{bNiceArray},
\end{aligned}
\end{equation}
where $\Sigma_E>0$ and $\widetilde{\Sigma}_R>0$; and let
\begin{equation*}
\begin{aligned}
\begin{bNiceArray}{c|c}[nullify-dots]
\widetilde{U}\widetilde{J}\widetilde{U}^T&\widetilde{U}\widetilde{B}
\end{bNiceArray}
=
\begin{bNiceArray}{ccc|ccc}[first-row,first-col,nullify-dots]
&\widetilde{n}_1&\widetilde{n}_2& \widetilde{n}_3&m\\
\widetilde{n}_1&\widetilde{J}_{11}&\widetilde{J}_{12}&\widetilde{J}_{13}&G_{1}\\
\widetilde{n}_2&-\widetilde{J}_{12}^T&\widetilde{J}_{22}&\widetilde{J}_{23}&G_{2}\\
\widetilde{n}_3&-\widetilde{J}_{13}^T&\widetilde{J}_{23}^T&\widetilde{J}_{33}&G_{3}\\
\end{bNiceArray}.
\end{aligned}
\end{equation*}

Applying \eqref{eq20241203150201} and \eqref{eq20241203160801}, we obtain
\begin{align*}
\rank\begin{bmatrix}
\widetilde{R}_{11}&\widetilde{R}_{12}&0&G_1\\
\widetilde{R}_{12}^T&\widetilde{\Sigma}_R&0&G_2\\
0&0&0&G_3\\
\end{bmatrix}=\rank\begin{bmatrix}
\widetilde{R}_{11}&\widetilde{R}_{12}\\
\widetilde{R}_{12}^T&\widetilde{\Sigma}_R\\
\end{bmatrix},
\end{align*}
and therefore $\rank(G_3)=0$, i.e., $G_3=0$.

Combining this fact with condition (1) and \eqref{eq20241203160801} yields
\begin{align*}
&\rank\begin{bmatrix}
\widetilde{E}&\widetilde{J}-(\widetilde{R}+\widetilde{B}K\widetilde{B}^T)\\
0&\widetilde{E}\\
0&\widetilde{B}^T\\
\end{bmatrix}=n+\rank(\widetilde{E})\\
&\Longrightarrow
\rank\begin{bmatrix}
\widetilde{J}_{22}-\widetilde{\Sigma}_R&\widetilde{J}_{23}\\
-\widetilde{J}_{23}^T&\widetilde{J}_{33}\\
G_2^T&0\\
\end{bmatrix}=n-\widetilde{n}_1\\
&\Longrightarrow
\rank\begin{bmatrix}
\widetilde{J}_{23}\\
\widetilde{J}_{33}\\
\end{bmatrix}=\widetilde{n}_3.
\end{align*}
Hence, by \eqref{eq20241203160801} and part (i) of Lemma \ref{leChu2024S},
$$
\begin{bmatrix}
\widetilde{J}_{22}-\widetilde{\Sigma}_R&\widetilde{J}_{23}\\
-\widetilde{J}_{23}^T&\widetilde{J}_{33}
\end{bmatrix}
$$
is nonsingular. 
Therefore,
\begin{align*}
\rank\begin{bmatrix}
\widetilde{E}&0\\
\widetilde{J}-(\widetilde{R}+\widetilde{B}K\widetilde{B}^T)&\widetilde{E}\\
\end{bmatrix}=n+\rank(\widetilde{E}).
\end{align*}
By Lemma \ref{thDai1989}, the matrix pair $(\widetilde{E},\widetilde{J}-(\widetilde{R}+\widetilde{B}K\widetilde{B}^T))$ is impulse-free. 

Together with the stability result above, this proves (a).~\hfill$~~\square$

For practical use, Theorem \ref{eq202507222239} leads to the following representation-free design procedure for pH descriptor systems.

\begin{algorithm}[H]
\caption{Representation-free proportional-gain design}
\label{ALPg}
\begin{algorithmic}[1]
\Require $E$, $A$, $B$, $C$ satisfying Assumption \ref{assu:pH}
\Ensure proportional output feedback $u=-Ky+v$
\Statex
\State Check conditions (i)--(ii) in Theorem \ref{eq202507222239}.
\State If both hold, choose any $K>0$.
\end{algorithmic}
\end{algorithm}

Thus, proportional-feedback synthesis reduces to solvability verification: once the stated rank conditions hold, any positive definite gain is admissible.

\section{Coefficient-based proportional and derivative feedback design with prescribed rank}\label{Se4}
Unlike proportional feedback, derivative feedback alters the leading matrix $E$ and thus requires a genuinely constructive design. 
This section develops such a procedure. 
We derive decompositions from $E,A,B,C$ and use them to compute a derivative gain $F$ that assigns the desired dynamical order.

By Proposition~\ref{re202510181952}, the solvability of the considered problems implies that $Q$ is nonsingular in every pH representation. 
Therefore, throughout this section we assume that, for any pH representation under consideration,  $Q$ is nonsingular; this assumption is used only in the analysis and not as computational input.

\begin{lemma}\rm\label{le202507221916}
Consider the descriptor system \eqref{eqI011} under Assumption \ref{assu:pH}, with nonsingular $Q$. 
Then there exist orthogonal matrices $U$, $V$, and $W$ such that
\begin{equation*}
\begin{aligned}
&\begin{bNiceArray}{c|c}[nullify-dots]
UEV & UBW
\end{bNiceArray}
=
\begin{bNiceArray}{ccccc|ccc}[first-row,first-col,nullify-dots]
&n_1&n_2& n_3& n_4&n_5&n_1&m_1& n_3\\
n_1&\widetilde{E}_{11}&E_{12}&0&0&0&B_{11}& 0&B_{13}\\
n_2&E_{21}&E_{22}&0&0&0&0&0&B_{23}\\
n_3&0&0&0&0&0&0& 0&B_{33}\\
n_4&0&0&0&0&0&0&0&0\\
n_5&0&0&0&0&0&0&0&0\\
\end{bNiceArray},
\end{aligned}
\end{equation*}
\begin{small}
\begin{equation*}
\begin{aligned}
\begin{bNiceArray}{c}[nullify-dots]
UAV \\
\hline
W^TCV
\end{bNiceArray}
=
\begin{bNiceArray}{lllll}[first-row,first-col,nullify-dots]
&n_1&n_2& n_3& n_4& n_5\\
n_1&A_{11}&A_{12}&A_{13}&A_{14}&A_{15}\\
n_2&A_{21}&A_{22}&A_{23}&A_{24}&A_{25}\\
n_3&A_{31}&A_{32}&A_{33}&A_{34}&A_{35}\\
n_4&A_{41}&A_{42}&A_{43}&A_{44}&0\\
n_5&A_{51}&A_{52}&A_{53}&0&0\\
\hline
n_1&C_{11}& 0&0&0&0\\
m_1&0& 0&0&0&0\\
n_3&C_{31}& C_{32}&C_{33}&0&0\\
\end{bNiceArray},
\end{aligned}
\end{equation*}
\end{small}
\vspace{-0.5cm}

\noindent
where $\rank\begin{bmatrix}
\widetilde{E}_{11}&E_{12}\\
E_{21}&E_{22}
\end{bmatrix}=n_1+n_2$, $\rank(B_{11})=\rank(C_{11})=n_1$, $\rank(E_{22})=n_2$, $\rank(B_{33})=\rank(C_{33})=n_3$, and $\rank(A_{44})=n_4$.
\end{lemma}
\textbf{Proof.~~}
A constructive proof, directly implementable as a numerical procedure, is provided in Appendix \ref{Appendix01}. 
The required orthogonal transformations are computed from the coefficient matrices, whereas the existence of a pH representation is used only to justify the resulting condensed block structure.~\hfill$~~\square$

If general nonsingular equivalence transformations are allowed, the condensed form in Lemma \ref{le202507221916} can be further simplified.

\begin{lemma}\rm\label{eq202507221956}
Consider the descriptor system \eqref{eqI011} under Assumption \ref{assu:pH}, with nonsingular $Q$. 
Then there exist nonsingular matrices $S$, $T$, and an orthogonal matrix $W$ such that
\vspace{-0.5cm}
\begin{equation*}
\begin{aligned}
&\begin{bNiceArray}{c|c}[nullify-dots]
SET & SBW
\end{bNiceArray}
=
\begin{bNiceArray}{ccccc|ccc}[first-row,first-col,nullify-dots]
&n_1&n_2& n_3& n_4&n_5&n_1&m_1& n_3\\
n_1&E_{11}&0&0&0&0&B_{11}& 0&0\\
n_2&0&E_{22}&0&0&0&0&0&0\\
n_3&0&0&0&0&0&0& 0&B_{33}\\
n_4&0&0&0&0&0&0&0&0\\
n_5&0&0&0&0&0&0&0&0\\
\end{bNiceArray},
\end{aligned}
\end{equation*}
\vspace{-1cm}
\begin{small}
\begin{equation*}
\begin{aligned}
\begin{bNiceArray}{c}[nullify-dots]
SAT\\
\hline
W^TCT
\end{bNiceArray}
=
\begin{bNiceArray}{lllll}[first-row,first-col,nullify-dots]
&n_1&n_2& n_3& n_4& n_5\\
n_1&A_{11}&A_{12}&A_{13}&0&A_{15}\\
n_2&A_{21}&A_{22}&A_{23}&0&A_{25}\\
n_3&A_{31}&A_{32}&A_{33}&0&A_{35}\\
n_4&0&0&0&A_{44}&0\\
n_5&A_{51}&A_{52}&A_{53}&0&0\\
\hline
n_1&C_{11}& 0&0&0&0\\
m_1&0& 0&0&0&0\\
n_3&0& 0&C_{33}&0&0\\
\end{bNiceArray},
\end{aligned}
\end{equation*}
\end{small}
\hspace{-1em}
where $\rank(E_{11})=\rank(B_{11})=\rank(C_{11})=n_1$, $\rank(E_{22})=n_2$, $\rank(B_{33})=\rank(C_{33})=n_3$, and $\rank(A_{44})=n_4$.
\end{lemma}
\textbf{Proof.~~}
Apply block Gaussian elimination to the matrices in Lemma \ref{le202507221916}.~\hfill$~~\square$

To characterize the dynamical order $\rank(E+BFC)$ in closed-loop systems that are regular, impulse-free, and asymptotically stable, the following result is needed.

\begin{corollary}\label{cor20250905204501}
\rm
Let the descriptor system \eqref{eqI011} be a pH descriptor system in the condensed form of Lemma \ref{eq202507221956}.
Then

(i)
\vspace{-0.5cm}
$$
\rank\begin{bmatrix}
E&B
\end{bmatrix}
=
\rank\begin{bmatrix}
E\\
C
\end{bmatrix}
=n_1+n_2+n_3.
$$
and
\begin{align*}
&\rank\begin{bmatrix}
E&B\\
\end{bmatrix}
+\rank\begin{bmatrix}
E\\
C\\
\end{bmatrix}
-\rank\begin{bmatrix}
E&B\\
C&0\\
\end{bmatrix}\\
=&\rank\begin{bmatrix}
E&B\\
\end{bmatrix}
-\rank(B)=n_2.
\end{align*}

(ii) $n_5=0$ if and only if
\begin{align}\label{eq20250916145201}
\rank\begin{bmatrix}
E&0&0\\
A&E&B\\
C&0&0\\
\end{bmatrix}=n+\rank\begin{bmatrix}
E&B
\end{bmatrix}.
\end{align}

(iii)
$$
\rank\begin{bmatrix}
A_{15}\\
A_{35}\\
\end{bmatrix}=n_5
$$
if and only if the matrix $T_a^TAS_b$ has full column rank, where $T_a=T_{\infty}(ES_{\infty}(C))$ and
\begin{align*}
S_b=S_{\infty}\left(\begin{bmatrix}
E\\
(T_{\infty}\left(\begin{bmatrix}
E&B
\end{bmatrix}\right))^TA\\
C\\
\end{bmatrix}\right).
\end{align*}
\end{corollary}
\textbf{Proof.~~}
The proofs of (i) and (ii) follow by direct calculation.

For (iii), note that
\begin{small}
\begin{align*}
S_{\infty}(C)=T\begin{bmatrix}
0&0&0\\
I&0&0\\
0&0&0\\
0&I&0\\
0&0&I\\
\end{bmatrix},
~
T_{\infty}\left(\begin{bmatrix}
E&B
\end{bmatrix}\right)=S^T\begin{bmatrix}
0&0\\
0&0\\
0&0\\
I&0\\
0&I\\
\end{bmatrix},
\end{align*}
\end{small}
so
\vspace{-0.5cm}
\begin{small}
\begin{align*}
T_a=S^T\begin{bmatrix}
I&0&0&0\\
0&0&0&0\\
0&I&0&0\\
0&0&I&0\\
0&0&0&I\\
\end{bmatrix},
~
S_b=T\begin{bmatrix}
0\\
0\\
0\\
0\\
I\\
\end{bmatrix}.
\end{align*}
\end{small}
Therefore, $\rank(T_a^TAS_b)=\rank\begin{bmatrix}
A_{15}\\
A_{35}\\
\end{bmatrix}$ by direct calculation.~\hfill$~~\square$

Using these condensed forms, we first characterize the assignable range of dynamical orders in terms of the given matrices.

\begin{theorem}\rm\label{th20240604}
Consider the descriptor system \eqref{eqI011} under Assumption \ref{assu:pH}, with nonsingular $Q$. 
Let $r$ be an integer. Then the following statements are equivalent.

\noindent
(a) There exists a matrix $F\in\mathbb{R}^{m\times m}$ such that
$$
\rank(E+BFC)=r,
$$
and such that the closed-loop system
\begin{align}\label{eq202507222103}
\left\{\begin{aligned}
(E+BFC)\dot{x}&=Ax+Bv,\\
y&=Cx,
\end{aligned}\right.
\end{align}
still admits a pH representation.

\noindent
(b) The integer $r$ satisfies
\begin{align}\label{eq202507222055}
\rank\begin{bmatrix}
E&B
\end{bmatrix}
-\rank(B)\leq r\leq \rank\begin{bmatrix}
E&B
\end{bmatrix}.
\end{align}
\end{theorem}

\textbf{Proof.~~}
``(a) $\Longrightarrow$ (b)''
This follows from Lemma \ref{th2.4}, Lemma \ref{eq202507221956}, and Corollary \ref{cor20250905204501}.

``(b) $\Longrightarrow$ (a)''
Choose any pH representation of system \eqref{eqI011}. 
This representation is used only in the proof and is not required by the construction.

By Lemma \ref{eq202507221956}, write
\begin{align*}
(S^T)^{-1}QT=\begin{bNiceArray}{ccccc}[first-row,first-col,nullify-dots]
&n_1&n_2& n_3& n_4& n_5\\
n_1&Q_{11}&Q_{12}&Q_{13}&Q_{14}&Q_{15}\\
n_2&Q_{21}&Q_{22}&Q_{23}&Q_{24}&Q_{25}\\
n_3& Q_{31}&Q_{32}&Q_{33}&Q_{34}&Q_{35}\\
n_4&Q_{41}&Q_{42}&Q_{43}&Q_{44}&Q_{45}\\
n_5&Q_{51}&Q_{52}&Q_{53}&Q_{54}&Q_{55}\\
\end{bNiceArray}.
\end{align*}
From $Q^TE=E^TQ\geq0$ and $C=B^TQ$, we obtain
$$
\begin{bmatrix}
Q_{13}&Q_{14}&Q_{15}\\
Q_{23}&Q_{24}&Q_{25}\\
\end{bmatrix}=0
$$
and
\begin{align*}
Q_{11}=(B_{11}^T)^{-1}C_{11},\quad Q_{12}=0,\quad Q_{21}=0,\quad Q_{31}=0,\\
Q_{32}=0,\quad Q_{34}=0,\quad Q_{35}=0,\quad Q_{33}=(B_{33}^T)^{-1}C_{33}.
\end{align*}
Hence
\begin{align*}
(S^T)^{-1}QT=\begin{bmatrix}
(B_{11}^T)^{-1}C_{11}&0&0&0&0\\
0&Q_{22}&0&0&0\\
0&0&(B_{33}^T)^{-1}C_{33}&0&0\\
Q_{41}&Q_{42}&Q_{43}&Q_{44}&Q_{45}\\
Q_{51}&Q_{52}&Q_{53}&Q_{54}&Q_{55}\\
\end{bmatrix}.
\end{align*}

Next, let
\begin{align}\label{eq20250909163601}
W^TFW=\begin{bNiceArray}{ccc}[first-row,first-col,nullify-dots]
&n_1&m_1& n_3\\
n_1& F_{11}&F_{12}&F_{13}\\
m_1&F_{21}&F_{22}&F_{23}\\
n_3&F_{31}&F_{32}&F_{33}\\
\end{bNiceArray}.
\end{align}
Then
\begin{equation}\label{eq20250909163602}
\begin{aligned}
&SET+SBWW^TFWW^TCT\\
=&\begin{bmatrix}
E_{11}+B_{11}F_{11}C_{11}&0&B_{11}F_{13}C_{33}&0&0\\
0&E_{22}&0&0&0\\
B_{33}F_{31}C_{11}&0&B_{33}F_{33}C_{33}&0&0\\
0&0&0&0&0\\
0&0&0&0&0\\
\end{bmatrix}
\end{aligned}
\end{equation}
and
\begin{align*}
&((S^T)^{-1}QT)^T(SET+SBWW^TFWW^TCT)\\
=&\begin{bmatrix}
C_{11}^TB_{11}^{-1}E_{11}+C_{11}^TF_{11}C_{11}&0&C_{11}^TF_{13}C_{33}&0&0\\
0&Q_{22}^TE_{22}&0&0&0\\
C_{33}^TF_{31}C_{11}&0&C_{33}^TF_{33}C_{33}&0&0\\
0&0&0&0&0\\
0&0&0&0&0\\
\end{bmatrix}.
\end{align*}
Consequently,
\begin{equation}\label{eq202510202215}
\begin{aligned}
&\rank(E+BFC)=\rank(Q^T(E+BFC))\\
=&\rank\begin{bmatrix}
C_{11}^TB_{11}^{-1}E_{11}+C_{11}^TF_{11}C_{11}&C_{11}^TF_{13}C_{33}\\
C_{33}^TF_{31}C_{11}&C_{33}^TF_{33}C_{33}\\
\end{bmatrix}\\
&+\rank(E_{22})
\end{aligned}
\end{equation}
and
\begin{equation}\label{eq202510202216}
\begin{aligned}
&Q^T(E+BFC)\geq0
\Longleftrightarrow\\
&\begin{bmatrix}
C_{11}^TB_{11}^{-1}E_{11}+C_{11}^TF_{11}C_{11}&C_{11}^TF_{13}C_{33}\\
C_{33}^TF_{31}C_{11}&C_{33}^TF_{33}C_{33}\\
\end{bmatrix}\geq 0.
\end{aligned}
\end{equation}

Let the integer $r$ satisfy \eqref{eq202507222055}. Choose
\begin{equation}\label{eq202511172053}
\begin{aligned}
&\begin{bmatrix}
F_{11}&F_{13}\\
F_{31}&F_{33}\\
\end{bmatrix}
=
\begin{bmatrix}
C_{11}^T&0\\
0&C_{33}^T\\
\end{bmatrix}^{-1}\\
&\times\left(-\begin{bmatrix}
C_{11}^TB_{11}^{-1}E_{11}&0\\
0&0\\
\end{bmatrix}+\mathcal{F}\right)
\begin{bmatrix}
C_{11}&0\\
0&C_{33}\\
\end{bmatrix}^{-1},
\end{aligned}
\end{equation}
and define $\mathcal{F} \in \mathbb{R}^{(n_1+n_3)\times (n_1+n_3)}$, $\mathcal{F}\geq 0$, and $\rank(\mathcal{F})=r-\rank(E_{22})$. Then, by direct calculation using \eqref{eq202510202215} and \eqref{eq202510202216}, we obtain
$$
\rank(E+BFC)=r
~\text{and}~
Q^T(E+BFC)\geq0.
$$
This proves the sufficiency.~\hfill$~~\square$

Building on Theorem \ref{th20240604}, we now derive a coefficient-based synthesis result for Problem \ref{pr2025072143601}.

\begin{theorem}\rm\label{th20250916152201}
Consider the descriptor system \eqref{eqI011} under Assumption \ref{assu:pH}. 
Let $r$ be a prescribed integer. Then the following statements are equivalent.

\noindent
(a) There exist matrices $K,F\in\mathbb{R}^{m\times m}$ such that the matrix pair $(E+BFC,A-BKC)$ is regular, impulse-free, and asymptotically stable, satisfies
$$
\rank(E+BFC)=r,
$$
and such that the closed-loop system \eqref{eq202411160142} still admits a pH representation.

\noindent
(b) The following coefficient-level conditions hold:

(i) The matrix $T_a^TAS_b$ has full column rank.

(ii) Condition \eqref{eq20241203143301} holds for all $s\in\mathbb{C}$ with $\Re(s)=0$.

(iii) The integer $r$ satisfies
$$
\rank\begin{bmatrix}
E&B
\end{bmatrix}
-\rank(B)\leq r\leq 
\rank\begin{bmatrix}
E&B
\end{bmatrix}
-\rank(T_a^TAS_b).
$$
Here, $T_a$ and $S_b$ are defined in part \textup{(iii)} of Corollary \ref{cor20250905204501}.
\end{theorem}
\textbf{Proof.~~}
We prove the equivalence of \textup{(a)} and \textup{(b)} by combining the coefficient-level rank conditions with Theorem \ref{eq202507222239}.

``(a) $\Longrightarrow$ (b)''
Let $K$ and $F$ be such that the matrix pair $(E+BFC,A-BKC)$ is regular, impulse-free, and asymptotically stable. Then
\begin{equation}\label{eq202511172122}
\begin{aligned}
&\rank\begin{bmatrix}s(E+BFC)-A&B\end{bmatrix}\\
=&
\rank\begin{bmatrix}s(E+BFC)-A\\
C\end{bmatrix}=n
\end{aligned}
\end{equation}
for all $s\in \mathbb{C}$ with $\Re(s)=0$, and
\begin{equation}\label{eq20250911152001}
\begin{aligned}
&\rank\begin{bmatrix}
E+BFC&A\\
0& E+BFC\\
0&C\\
\end{bmatrix}\\
=&
n+\rank(E+BFC),
\end{aligned}
\end{equation}
by Theorem \ref{eq202507222239}. Clearly, condition (ii) follows from \eqref{eq202511172122}.

We now prove condition (i). 
By Proposition \ref{re202510181952}, we have $\rank(Q)=n$. 
Part (iii) of Corollary \ref{cor20250905204501} shows that it suffices to prove
\begin{align}\label{eq202511172135}
\rank\begin{bmatrix}
A_{15}\\
A_{35}\\
\end{bmatrix}=n_5.
\end{align}
Using Lemma \ref{eq202507221956}, \eqref{eq20250909163601}, and \eqref{eq20250909163602}, we obtain
\begin{equation}\label{eq20250909163901}
\begin{aligned}
&\rank\begin{bmatrix}
E+BFC&A\\
0& E+BFC\\
0&C\\
\end{bmatrix}\\
=&\rank\begin{bmatrix}
E_{11}+B_{11}F_{11}C_{11}&B_{11}F_{13}C_{33}&A_{15}\\
B_{33}F_{31}C_{11}&B_{33}F_{33}C_{33}&A_{35}\\
\end{bmatrix}\\
&+n_1+n_2+n_3+n_4+n_2
\end{aligned}
\end{equation}
and
\begin{equation}\label{eq20251020163901}
\begin{aligned}
&\rank(E+BFC)\\
=&
n_2+\rank\begin{bmatrix}
E_{11}+B_{11}F_{11}C_{11}&B_{11}F_{13}C_{33}\\
B_{33}F_{31}C_{11}&B_{33}F_{33}C_{33}\\
\end{bmatrix}.
\end{aligned}
\end{equation}
Substituting \eqref{eq20250909163901} and \eqref{eq20251020163901} into \eqref{eq20250911152001} yields
\begin{equation}\label{202511172232}
\begin{aligned}
&\rank\begin{bmatrix}
E_{11}+B_{11}F_{11}C_{11}&B_{11}F_{13}C_{33}&A_{15}\\
B_{33}F_{31}C_{11}&B_{33}F_{33}C_{33}&A_{35}\\
\end{bmatrix}\\
=&\rank\begin{bmatrix}
E_{11}+B_{11}F_{11}C_{11}&B_{11}F_{13}C_{33}\\
B_{33}F_{31}C_{11}&B_{33}F_{33}C_{33}\\
\end{bmatrix}+n_5,
\end{aligned}
\end{equation}
and therefore \eqref{eq202511172135} holds.

By \eqref{eq20250909163901} and Lemma \ref{th2.4}, we get
\begin{align*}
&n_1+n_2+n_3+n_4+n_2+\rank\begin{bmatrix}
A_{15}\\
A_{35}\\
\end{bmatrix}\\
\leq& \rank\begin{bmatrix}
E+BFC&A\\
0& E+BFC\\
0&C\\
\end{bmatrix}\\
\leq& n_1+n_2+n_3+n_4+n_2+n_1+n_3.
\end{align*}
Combining this with \eqref{eq20250911152001}, we obtain condition (iii).

``(b) $\Longrightarrow$ (a)''
Assume that conditions (i)-(iii) hold. 
We prove sufficiency by constructing $F$ from the coefficient data and then applying Theorem \ref{eq202507222239} to choose $K>0$.

For the proof, fix an arbitrary pH representation; condition (ii) ensures that its $Q$ is nonsingular.

By Theorem \ref{eq202507222239}, it remains to show that there exists $F$ such that \eqref{eq202511172122}, \eqref{eq20250911152001}, $\rank(E+BFC)=r$, and $Q^T(E+BFC)\geq 0$ hold for any integer $r$ satisfying condition (iii). 
For the feedback gain $K$, we may choose any positive definite matrix.

Since \eqref{eq202511172122} holds for any matrix $F$ satisfying the rank and positive-semidefiniteness constraints, the problem reduces to constructing an $F$ that also satisfies \eqref{eq20250911152001}. 
Furthermore, when condition (ii) holds, Theorem \ref{eq202507222239} guarantees that $Q$ is nonsingular. 
Consequently, by the proof of necessity, \eqref{eq20250911152001} is equivalent to \eqref{202511172232}. 
Therefore, it remains to construct a matrix
$$
\begin{bmatrix}
F_{11}&F_{13}\\
F_{31}&F_{33}
\end{bmatrix}
$$
with the required properties.

By Lemma \ref{eq202507221956}, \eqref{eq20250909163601}, and \eqref{eq20250909163602}, we have
\begin{small}
\begin{flalign*}
&\hspace{-4mm}\rank\begin{bmatrix}
E_{11}+B_{11}F_{11}C_{11}&B_{11}F_{13}C_{33}&A_{15}\\
B_{33}F_{31}C_{11}&B_{33}F_{33}C_{33}&A_{35}\\
\end{bmatrix}=\\
&\hspace{-4mm}\rank \left(\begin{bmatrix}
C_{11}^TB_{11}^{-1} &0\\
0&C_{33}^TB_{33}^{-1}\\
\end{bmatrix}
\begin{bNiceArray}{ccc}[nullify-dots]
E_{11}+B_{11}F_{11}C_{11}&B_{11}F_{13}C_{33}&A_{15}\\
B_{33}F_{31}C_{11}&B_{33}F_{33}C_{33}&A_{35}\\
\end{bNiceArray}\right)\\
&\hspace{-4mm}=\rank\begin{bmatrix}
C_{11}^TB_{11}^{-1}E_{11}+C_{11}^TF_{11}C_{11}&C_{11}^TF_{13}C_{33}&C_{11}^TB_{11}^{-1}A_{15}\\
C_{33}^TF_{31}C_{11}&C_{33}^TF_{33}C_{33}&C_{33}^TB_{33}^{-1}A_{35}\\
\end{bmatrix}.
\end{flalign*}
\end{small}
From condition (i), there exists an orthogonal matrix $\widetilde{U}$ such that
\begin{align*}
\widetilde{U}\begin{bmatrix}
C_{11}^TB_{11}^{-1}A_{15}\\
C_{33}^TB_{33}^{-1}A_{35}\\
\end{bmatrix}
=
\begin{bNiceArray}{c}[first-row,first-col,nullify-dots]
&n_5\\
n_5&\widetilde{\Sigma_A}\\
\widetilde{n}_5&0\\
\end{bNiceArray},
\end{align*}
where $\rank(\widetilde{\Sigma_A})=n_5$ and $\widetilde{n}_5=\rank(B)-n_5$.

For any integer $r$ satisfying condition (iii), we may choose $\begin{bmatrix}
F_{11}&F_{13}\\
F_{31}&F_{33}\\
\end{bmatrix}$
such that
\begin{equation}\label{20260422211601}
\begin{aligned}
&\widetilde{U}\begin{bmatrix}
C_{11}^T&0\\
0&C_{33}^T\\
\end{bmatrix}
\begin{bmatrix}
F_{11}&F_{13}\\
F_{31}&F_{33}\\
\end{bmatrix}
\begin{bmatrix}
C_{11}&0\\
0&C_{33}\\
\end{bmatrix}\widetilde{U}^T\\
=&-\widetilde{U}\begin{bmatrix}
C_{11}^TB_{11}^{-1}E_{11}&0\\
0&0\\
\end{bmatrix}\widetilde{U}^T
+\begin{bNiceArray}{cc}[first-row,first-col,nullify-dots]
&n_5&\widetilde{n}_5\\
n_5&0&0\\
\widetilde{n}_5&0&\mathcal{F}_{33}\\
\end{bNiceArray},
\end{aligned}
\end{equation}
where $\mathcal{F}_{33}\geq 0$ and $\rank(\mathcal{F}_{33})=r-\rank(E_{22})$. A direct calculation then shows that \eqref{202511172232}, $\rank(E+BFC)=r$, and $Q^T(E+BFC)\geq 0$ all hold.~\hfill$~~\square$

The following corollary characterizes when every admissible rank in \eqref{eq202507222055} can be assigned.

\begin{corollary}\rm
Consider the descriptor system \eqref{eqI011} under Assumption \ref{assu:pH}. 
Then every integer $r$ satisfying \eqref{eq202507222055} is assignable by some $K,F\in\mathbb{R}^{m\times m}$ such that $(E+BFC,A-BKC)$ is regular, impulse-free, and asymptotically stable, and the closed-loop system \eqref{eq202411160142} still admits a pH representation, if and only if \eqref{eq20250916145201} holds and \eqref{eq20241203143301} holds for all $s\in \mathbb{C}$ with $\Re(s)=0$.
\end{corollary}
\textbf{Proof.~~}By Proposition \ref{re202510181952} and Theorem \ref{th20250916152201}, the matrix $Q$ in any pH representation of \eqref{eqI011} is nonsingular under the hypotheses of this corollary.

``$\Longrightarrow$'' Let $\rank(E+BFC)=\rank\begin{bmatrix}
E&B
\end{bmatrix}$.
Then $n_5=0$ follows from \eqref{202511172232}. The necessity follows from part (ii) of Corollary \ref{cor20250905204501} and Theorem \ref{th20250916152201}.

``$\Longleftarrow$''
Because \eqref{eq20250916145201} holds, we also have $n_5=0$. Hence \eqref{202511172232} holds automatically for any $F$ such that $\rank(E+BFC)=r$ and $Q^T(E+BFC)\geq0$, where $r$ is any integer satisfying \eqref{eq202507222055}. Therefore, the sufficiency follows from Theorem \ref{th20250916152201}.~\hfill$~~\square$

For a prescribed integer $r$, Theorem \ref{th20250916152201} yields the following coefficient-based design procedure.

\begin{algorithm}[H]
\caption{Coefficient-based  proportional-derivative feedback design with prescribed rank}
\label{AlPDg}
\begin{algorithmic}[1]
\Require $E$, $A$, $B$, $C$ satisfying Assumption \ref{assu:pH}
\Ensure proportional-derivative output feedback law $u=-Ky-F\dot{y}+v$
\Statex
\State Verify \eqref{eq20241203143301}.
\State Compute the condensed form in Lemma \ref{le202507221916} and Lemma \ref{eq202507221956}.
\State Compute $T_a$ and $S_b$ in part \textup{(iii)} of Corollary \ref{cor20250905204501}.
\State Verify that $T_a^TAS_b$ has full column rank.
\State Check the feasibility conditions in Theorem \ref{th20250916152201} for the prescribed rank $r$.
\State If these conditions hold, construct the block matrix according to \eqref{20260422211601}, and choose the remaining blocks of the derivative gain $F$ arbitrarily with compatible dimensions.
\State Recover the derivative feedback matrix $F$ via the transformation formula \eqref{eq20250909163601}.
\State Choose any $K>0$.
\end{algorithmic}
\end{algorithm}

Thus, the derivative gain $F$ assigns the prescribed dynamical order, while $K$ can be any positive definite matrix once solvability is verified.
The rank conditions involving $s$ can be checked numerically via generalized eigenvalue methods (e.g., QZ algorithm). 
The orthogonal transformations in Algorithm~\ref{AlPDg} are standard numerical linear algebra procedures.

The following example illustrates the direct verification of the coefficient-level conditions and the construction of feedback gains without a pH representation.

\begin{example}\rm
Consider the descriptor system \eqref{eqI011} with
$$
E=\begin{bmatrix}1&0\\0&0\end{bmatrix},~
A=\begin{bmatrix}0&1\\-1&0\end{bmatrix},~
B=\begin{bmatrix}1\\1\end{bmatrix},~
C=\begin{bmatrix}2&1\end{bmatrix}.
$$

The system is completely controllable, completely observable, and positive real. 
Therefore, by the realization result in \citep{Chu2025A}, it admits at least one pH representation, although no such representation is constructed here.

Although the open-loop system is neither impulse-free nor asymptotically stable, the coefficient-level conditions of Theorem \ref{th20250916152201} are satisfied for $r\in\{1,2\}$. 
Hence admissible feedback matrices can be selected directly from the coefficient data. 

The following feedback choices are admissible.

For $r=1$, Algorithm \ref{AlPDg} allows the choice $F=0$ and $K=1$. 
The resulting closed-loop system is given by
$$
E=\begin{bmatrix}1&0\\0&0\end{bmatrix},~
A-BKC=\begin{bmatrix}-2&0\\-3&-1\end{bmatrix},~
B=\begin{bmatrix}1\\1\end{bmatrix},~
C=\begin{bmatrix}2&1\end{bmatrix}.
$$

For $r=2$, Algorithm \ref{AlPDg} allows the choice $K=1$ and $F=1$. 
Then the closed-loop system becomes
$$
E+BFC=\begin{bmatrix}3&1\\2&1\end{bmatrix},~
A-BKC=\begin{bmatrix}-2&0\\-3&-1\end{bmatrix},~
B=\begin{bmatrix}1\\1\end{bmatrix},
$$
and $C=\begin{bmatrix}2&1\end{bmatrix}$.

In both cases, the closed-loop system is regular, impulse-free, asymptotically stable, and port-Hamiltonian. 

This example shows that the coefficient-level conditions of Theorem~\ref{th20250916152201} allow stabilizing output feedback design without first computing a pH representation. 
By contrast, applying the results of \citep{Chu2024S} would require an explicit representation.
\end{example}

\section{Conclusion}
This paper has presented a coefficient-level framework for output-feedback stabilization of pH descriptor systems under the assumption that a pH representation exists. 
For proportional feedback, solvability is characterized by rank conditions involving only the coefficient matrices, and any symmetric positive definite gain is admissible. 
For derivative feedback, a coefficient-based construction assigns the prescribed dynamical order without forming a pH representation. 
Future work will focus on efficient decomposition algorithms and verifiable conditions for the existence of pH representations.
\begin{ack}
The authors are grateful to the anonymous reviewers for their careful reading and constructive comments, which have significantly improved the presentation of this paper.

This work was supported by the National Key Research and Development Program of China (Grant No. 2023YFB3001604).
\end{ack}

\bibliographystyle{plainnat}
\bibliography{autosam}

\appendix
\section{Proof of Lemma \ref{le202507221916}}\label{Appendix01}
This appendix provides a constructive proof that is numerically implementable. 
We choose an arbitrary pH representation of system \eqref{eqI011}, using it only to justify the block structure.

Let $r_e=\rank(E)$ and $r_b=\rank(B)$. 
Since $Q$ is nonsingular and $Q^TP=0$, we obtain $P=0$ for the pH descriptor system \eqref{eqI011}. 
Moreover, since $C=B^TQ$, we have $\rank(C)=\rank(B)=r_b$. 
We construct orthogonal matrices $U$, $V$, and $W$ by using QR decompositions and singular value decompositions.

\textit{Step 1.}
Compute orthogonal matrices $U_1$ and $V_1$ such that
\begin{equation*}
E^{(1)}=U_1EV_1=\begin{bNiceArray}{cc}[first-row,first-col,nullify-dots]
&r_{e}&\widetilde{n}\\
r_{e}&\widetilde{\Sigma_E}&0\\
\widetilde{n}&0&0\\
\end{bNiceArray},
\end{equation*}
where $\widetilde{n}=n-r_e$ and $\rank(\widetilde{\Sigma_E})=r_e$.
Define
$$
B^{(1)}=U_1B=\begin{bNiceArray}{c}[first-row,first-col,nullify-dots]
&m\\
r_{e}&B_{11}^{(1)}\\
\widetilde{n}&B_{21}^{(1)}\\
\end{bNiceArray},~
C^{(1)}=CV_1=\begin{bNiceArray}{cc}[first-row,first-col,nullify-dots]
&r_{e}&\widetilde{n}\\
m&C_{11}^{(1)}& C_{21}^{(1)}\\
\end{bNiceArray}.
$$
From $Q^TE=E^TQ\geq0$ and $\rank(Q)=n$, it follows that
\begin{equation*}
Q^{(1)}=U_1QV_1=\begin{bNiceArray}{cc}[first-row,first-col,nullify-dots]
&r_{e}&\widetilde{n}\\
r_{e}&Q_{11}^{(1)}&0\\
\widetilde{n}&Q_{21}^{(1)}&Q_{22}^{(1)}\\
\end{bNiceArray},
\end{equation*}
where $\rank(Q_{11}^{(1)})=r_{e}$ and $\rank(Q_{22}^{(1)})=\widetilde{n}$.

\textit{Step 2.}
Compute orthogonal matrices $U_2$ and $W_1$ such that
\begin{equation*}
U_2B_{21}^{(1)}W_1=\begin{bNiceArray}{cc}[first-row,first-col,nullify-dots]
&m-n_3& n_3\\
n_3&0& B_{33}\\
\widetilde{n}_3&0&0\\
\end{bNiceArray},
\end{equation*}
where $\widetilde{n}_3=n-r_e-n_3$ and $\rank(B_{33})=n_3$.
Define
\begin{equation*}
B^{(2)}=\begin{bmatrix}
I&0\\
0&U_2\\
\end{bmatrix}B^{(1)}W_1,~
Q^{(2)}=\begin{bmatrix}
I&0\\
0&U_2\\
\end{bmatrix}Q^{(1)},
\end{equation*}
and $C^{(2)}=W_1^TC^{(1)}$.
That is,
\begin{equation*}
B^{(2)}=\begin{bmatrix}
B_{11}^{(1)}W_1\\
U_2B_{21}^{(1)}W_1\\
\end{bmatrix}
=
\begin{bNiceArray}{cc}[first-row,first-col,nullify-dots]
&m-n_3& n_3\\
r_{e}&B_{11}^{(2)}&B_{21}^{(2)}\\
n_3&0& B_{33}\\
\widetilde{n}_3&0&0\\
\end{bNiceArray},
\end{equation*}
where $\rank(B_{33})=n_3$;
\begin{equation*}
		Q^{(2)}=\begin{bmatrix}
			Q_{11}^{(1)}&0\\
			U_2Q_{21}^{(1)}&U_2Q_{22}^{(1)}\\
		\end{bmatrix}=\begin{bNiceArray}{ccccc|ccc}[first-row,first-col,nullify-dots]&r_{e}
                 &n_3&\widetilde{n}_3\\
			r_{e}&Q_{11}^{(2)}&0&0\\
			n_3&Q_{21}^{(2)}&Q_{22}^{(2)}&Q_{23}^{(2)}\\
			\widetilde{n}_3&Q_{31}^{(2)}&Q_{32}^{(2)}&Q_{33}^{(2)}\\
		\end{bNiceArray},
\end{equation*}
where $\rank(Q_{11}^{(2)})=r_{e}$ and $\rank\begin{bmatrix}
	Q_{22}^{(2)}&Q_{23}^{(2)}\\
	Q_{32}^{(2)}&Q_{33}^{(2)}\\
\end{bmatrix}=\widetilde{n}$;
\begin{equation*}
		C^{(2)}\hspace{-1mm}=\hspace{-1mm}\begin{bmatrix}
			W_1^TC_{11}^{(1)}& W_1^TC_{21}^{(1)}
		\end{bmatrix}\hspace{-1mm}	
=\begin{bNiceArray}{ccccc|ccc}[first-row,first-col,nullify-dots]&r_{e}&n_3&\widetilde{n}_3\\
			m-n_3&C_{11}^{(2)}& C_{12}^{(2)}&C_{13}^{(2)}\\
			n_3&C_{21}^{(2)}& C_{22}^{(2)}&C_{23}^{(2)}\\
		\end{bNiceArray}.
\end{equation*}
Since $C=B^TQ$, we further obtain
\begin{align*}
C^{(2)}=
\begin{bNiceArray}{ccc}[first-row,first-col,nullify-dots]
&r_{e}&n_3&\widetilde{n}_3\\
m-n_3&C_{11}^{(2)}& 0&0\\
n_3&C_{21}^{(2)}& C_{22}^{(2)}&C_{23}^{(2)}\\
\end{bNiceArray},
\end{align*}
where
$$
\rank\begin{bmatrix}C_{22}^{(2)}&C_{23}^{(2)}\end{bmatrix}
=
\rank\begin{bmatrix}(B_{33})^TQ_{22}^{(2)}&(B_{33})^TQ_{23}^{(2)}\end{bmatrix}
=n_3.
$$

\textit{Step 3.}
Compute an orthogonal matrix $V_2$ such that
\begin{equation*}
\begin{bmatrix}
C_{22}^{(2)}&C_{23}^{(2)}
\end{bmatrix}V_2=
\begin{bNiceArray}{cc}[first-row,first-col,nullify-dots]
&n_3&\widetilde{n}_3\\
n_3&C_{33} &0\\
\end{bNiceArray},
\end{equation*}
where $\rank(C_{33})=n_3$.
Define
\begin{equation*}
C^{(3)}=C^{(2)}\begin{bmatrix}
I&0\\
0&V_2\\
\end{bmatrix},~
E^{(3)}=\begin{bmatrix}
I&0\\
0&U_2\\
\end{bmatrix}E^{(1)}\begin{bmatrix}
I&0\\
0&V_2\\
\end{bmatrix},
\end{equation*}
$B^{(3)}=B^{(2)}$, and set
$$
Q^{(3)}=Q^{(2)}\begin{bmatrix}
I&0\\
0&V_2\\
\end{bmatrix}.
$$
Then
\begin{equation*}
		C^{(3)}=\begin{bmatrix}
			C_{11}^{(2)}&\begin{bmatrix}0&0\end{bmatrix}\\
			C_{21}^{(2)}&\begin{bmatrix}C_{22}^{(2)}&C_{23}^{(2)}\end{bmatrix}V_2\\
		\end{bmatrix}=\begin{bNiceArray}{ccccc|ccc}[first-row,first-col,nullify-dots]
                 &r_{e}&n_3&\widetilde{n}_3\\
			m-n_3&C_{11}^{(3)}& 0&0\\
			n_3&C_{21}^{(3)}& C_{33}&0\\
		\end{bNiceArray},
\end{equation*}
where $\rank(C_{33})=n_3$ and $\rank(C_{11}^{(3)})=r_b-n_3$;
\begin{equation*}
B^{(3)}=
\begin{bNiceArray}{cc}[first-row,first-col,nullify-dots]
&m-n_3& n_3\\
r_{e}&B_{11}^{(3)}&B_{21}^{(3)}\\
n_3&0& B_{33}\\
\widetilde{n}_3&0&0\\
\end{bNiceArray},
\end{equation*}
where $\rank(B_{33})=n_3$ and $\rank(B_{11}^{(3)})=r_b-n_3$;
\begin{equation*}
		E^{(3)}=\begin{bmatrix}
			\widetilde{\Sigma_E}&0\\
			0&0\\
		\end{bmatrix}=\begin{bNiceArray}{ccc}[first-row,first-col,nullify-dots]
                 &r_{e}&n_3&\widetilde{n}_3\\
			r_{e}&\widetilde{\Sigma_E}&0&0\\
			n_3&0&0&0\\
			\widetilde{n}_3&0&0&0\\
		\end{bNiceArray},
\end{equation*}
where $\rank(\widetilde{\Sigma_E})=r_e$;
\begin{equation*}
		Q^{(3)}=\begin{bmatrix}
			Q_{11}^{(2)}&{\begin{bmatrix}0&0\end{bmatrix}}\\{\begin{bmatrix}Q_{21}^{(2)}\\Q_{31}^{(2)}\end{bmatrix} }& {\begin{bmatrix}Q_{22}^{(2)}&Q_{23}^{(2)}\\
					Q_{32}^{(2)}&Q_{33}^{(2)}\end{bmatrix}}V_2\\
		\end{bmatrix}
=\begin{bNiceArray}{ccc}[first-row,first-col,nullify-dots]
                 &r_{e}&n_3&\widetilde{n}_3\\
			r_{e}&Q_{11}^{(3)}&0&0\\
			n_3&Q_{21}^{(3)}&Q_{22}^{(3)}&Q_{23}^{(3)}\\
			\widetilde{n}_3&Q_{31}^{(3)}&Q_{32}^{(3)}&Q_{33}^{(3)}\\
		\end{bNiceArray},
\end{equation*}
where $\rank(Q_{11}^{(3)})=r_{e}$ and
$$
\rank\begin{bmatrix}
Q_{22}^{(3)}&Q_{23}^{(3)}\\
Q_{32}^{(3)}&Q_{33}^{(3)}
\end{bmatrix}=\widetilde{n}.
$$
Because $C=B^TQ$, we obtain $Q_{23}^{(3)}=0$. Hence
\begin{equation*}
Q^{(3)}=
\begin{bNiceArray}{ccc}[first-row,first-col,nullify-dots]
&r_{e}&n_3&\widetilde{n}_3\\
r_{e}&Q_{11}^{(3)}&0&0\\
n_3&Q_{21}^{(3)}&Q_{22}^{(3)}&0\\
\widetilde{n}_3&Q_{31}^{(3)}&Q_{32}^{(3)}&Q_{33}^{(3)}\\
\end{bNiceArray},
\end{equation*}
where $\rank(Q_{11}^{(3)})=r_{e}$, $\rank(Q_{22}^{(3)})=n_3$, and $\rank(Q_{33}^{(3)})=\widetilde{n}_3$.

\textit{Step 4.}
Compute orthogonal matrices $U_3$ and $W_2$ such that
\begin{equation*}
U_3B_{11}^{(3)}W_2=
\begin{bNiceArray}{cc}[first-row,first-col,nullify-dots]
&n_1&m_1\\
n_1&B_{11}& 0\\
n_2&0&0\\
\end{bNiceArray},
\end{equation*}
where $\rank(B_{11})=r_b-n_3=n_1$, $n_2=r_e-n_1$, and $m_1=m-n_1-n_3$.
Define
$$
\widetilde{U}=\begin{bmatrix}
U_3&0&0\\
0&I&0\\
0&0&I\\
\end{bmatrix}
\begin{bmatrix}
I&0\\
0&U_2\\
\end{bmatrix}U_1,
~
\widetilde{W}=W_1\begin{bmatrix}
W_2&0\\
0&I
\end{bmatrix}.
$$
Then
\begin{equation*}
		\widetilde{U}B\widetilde{W}=\begin{bmatrix}
			U_3B_{11}^{(3)}W_2&U_3B_{21}^{(3)}\\
			0& B_{33}\\
			0&0\\
		\end{bmatrix}=\begin{bNiceArray}{ccc}[first-row,first-col,nullify-dots]
               &n_1&m_1& n_3\\
			n_1&B_{11}& 0&B_{13}\\
			n_2&0&0&B_{23}\\
			n_3&0& 0&B_{33}\\
			\widetilde{n}_4&0&0&0\\
		\end{bNiceArray},
\end{equation*}
where $\rank(B_{11})=n_1$, $\rank(B_{33})=n_3$, and $\widetilde{n}_4=n-n_1-n_2-n_3$.
Define
$$
Q^{(4)}=\begin{bmatrix}
U_3&0&0\\
0&I&0\\
0&0&I\\
\end{bmatrix}Q^{(3)},
~
C^{(4)}=\begin{bmatrix}
W_2^T&0\\
0&I
\end{bmatrix}C^{(3)}.
$$
Then
\begin{align*}
		Q^{(4)}&=\begin{bmatrix}
			U_3Q_{11}^{(3)}&0&0\\
			Q_{21}^{(3)}&Q_{22}^{(3)}&0\\
			Q_{31}^{(3)}&Q_{32}^{(3)}&Q_{33}^{(3)}\\
		\end{bmatrix}\\
&=\begin{bNiceArray}{cccc}[first-row,first-col,nullify-dots]
               &n_1&n_2& n_3& \widetilde{n}_4\\
			n_1&Q_{11}^{(4)}&Q_{12}^{(4)}&0&0\\
			n_2&Q_{21}^{(4)}&Q_{22}^{(4)}&0&0\\
			n_3&Q_{31}^{(4)}&Q_{32}^{(4)}&Q_{33}^{(4)}&0\\
			\widetilde{n}_4&Q_{41}^{(4)}&Q_{42}^{(4)}&Q_{43}^{(4)}&Q_{44}^{(4)}\\
		\end{bNiceArray},
\end{align*}
where
$\rank\begin{bmatrix}
	Q_{11}^{(4)}&Q_{12}^{(4)}\\
	Q_{21}^{(4)}&Q_{22}^{(4)}\\
\end{bmatrix}=n_1+n_2,~\rank(Q_{33}^{(4)})=n_3$, and $\rank(Q_{44}^{(4)})=\widetilde{n}_4$;
\begin{equation*}
		C^{(4)}=\begin{bmatrix}
			W_2^TC_{11}^{(3)}& 0&0\\
			C_{21}^{(3)}& C_{33}&0\\
		\end{bmatrix}=\begin{bNiceArray}{cccc}[first-row,first-col,nullify-dots]
               &n_1&n_2& n_3& \widetilde{n}_4\\
			n_1&C_{11}^{(4)}& C_{12}^{(4)}&0&0\\
			m_1&C_{21}^{(4)}& C_{22}^{(4)}&0&0\\
			n_3&C_{31}^{(4)}& C_{32}^{(4)}&C_{33}&0\\
		\end{bNiceArray}.
\end{equation*}
Since $C=B^TQ$, it follows that
\begin{equation*}
C^{(4)}=
\begin{bNiceArray}{cccc}[first-row,first-col,nullify-dots]
&n_1&n_2& n_3& \widetilde{n}_4\\
n_1&C_{11}^{(4)}& C_{12}^{(4)}&0&0\\
m_1&0& 0&0&0\\
n_3&C_{31}^{(4)}& C_{32}^{(4)}&C_{33}&0\\
\end{bNiceArray},
\end{equation*}
where
$$
\rank\begin{bmatrix}C_{11}^{(4)}& C_{12}^{(4)}\end{bmatrix}
=
\rank\begin{bmatrix}(B_{11})^TQ_{11}^{(4)}&(B_{11})^TQ_{12}^{(4)}\end{bmatrix}
=n_1.
$$

\textit{Step 5.}
Compute an orthogonal matrix $V_3$ such that
\begin{equation*}
\begin{bmatrix}
C_{11}^{(4)}& C_{12}^{(4)}
\end{bmatrix}V_3=
\begin{bNiceArray}{cc}[first-row,first-col,nullify-dots]
&n_1&n_2\\
n_1&C_{11}&0\\
\end{bNiceArray},
\end{equation*}
where $\rank(C_{11})=n_1$.
Define
$$
\widetilde{V}=V_1
\begin{bmatrix}
I&0\\
0&V_2\\
\end{bmatrix}
\begin{bmatrix}
V_3&0&0\\
0&I&0\\
0&0&I\\
\end{bmatrix}.
$$
Then
\begin{align*}
		\widetilde{W}^TC\widetilde{V}&=\begin{bmatrix}
			\begin{bmatrix}C_{11}^{(4)}& C_{12}^{(4)}\end{bmatrix}V_3&0&0\\
			\begin{bmatrix}0& 0\end{bmatrix}&0&0\\
			\begin{bmatrix}C_{31}^{(4)}& C_{32}^{(4)}\end{bmatrix}V_3&C_{33}&0
		\end{bmatrix}\\
&=\begin{bNiceArray}{cccc}[first-row,first-col,nullify-dots]
               &n_1&n_2& n_3& \widetilde{n}_4\\
			n_1&C_{11}& 0&0&0\\
			m_1&0& 0&0&0\\
			n_3&C_{31}& C_{32}&C_{33}&0\\
		\end{bNiceArray},
\end{align*}
where $\rank(C_{11})=n_1$, $\rank(C_{33})=n_3$, and
\begin{align*}
\widetilde{U}Q\widetilde{V}&=\begin{bmatrix}
			\begin{bmatrix}Q_{11}^{(4)}&Q_{12}^{(4)}\\
				Q_{21}^{(4)}&Q_{22}^{(4)}
			\end{bmatrix}V_3&\begin{bmatrix}0\\0\end{bmatrix}&\begin{bmatrix}0\\0\end{bmatrix}\\
			~\\
			\begin{bmatrix}Q_{31}^{(4)}&Q_{32}^{(4)}\end{bmatrix}V_3&Q_{33}^{(4)}&0\\
			~\\
			\begin{bmatrix}Q_{41}^{(4)}&Q_{42}^{(4)}\end{bmatrix}V_3&Q_{43}^{(4)}&Q_{44}^{(4)}\\
		\end{bmatrix}\\
&=\begin{bNiceArray}{cccc}[first-row,first-col,nullify-dots]
               &n_1&n_2& n_3& \widetilde{n}_4\\
			n_1&Q_{11}&Q_{12}&0&0\\
			n_2&Q_{21}&Q_{22}&0&0\\
			n_3&Q_{31}&Q_{32}&Q_{33}&0\\
			\widetilde{n}_4&Q_{41}&Q_{42}&Q_{43}&Q_{44}\\
		\end{bNiceArray},
\end{align*}
where $\rank\begin{bmatrix}
	Q_{11}&Q_{12}\\
	Q_{21}&Q_{22}\\
\end{bmatrix}=n_1+n_2$, $\rank(Q_{33})=n_3$, and $\rank(Q_{44})=\widetilde{n}_4$.

Since $C=B^TQ$, we obtain $Q_{12}=0$. Hence
\begin{equation*}
\widetilde{U}Q\widetilde{V}=
\begin{bNiceArray}{cccc}[first-row,first-col,nullify-dots]
&n_1&n_2& n_3& \widetilde{n}_4\\
n_1&Q_{11}&0&0&0\\
n_2&Q_{21}&Q_{22}&0&0\\
n_3&Q_{31}&Q_{32}&Q_{33}&0\\
\widetilde{n}_4&Q_{41}&Q_{42}&Q_{43}&Q_{44}\\
\end{bNiceArray},
\end{equation*}
where $\rank(Q_{11})=n_1$, $\rank(Q_{22})=n_2$, $\rank(Q_{33})=n_3$, and $\rank(Q_{44})=\widetilde{n}_4$.
Furthermore,
\begin{equation*}
		\widetilde{U}E\widetilde{V}=\begin{bmatrix}
			U_3\widetilde{\Sigma_E}V_3&0&0\\
			0&0&0\\
			0&0&0\\
		\end{bmatrix}=\begin{bNiceArray}{cccc}[first-row,first-col,nullify-dots]
               &n_1&n_2& n_3& \widetilde{n}_4\\
			n_1&\widetilde{E}_{11}&E_{12}&0&0\\
			n_2&E_{21}&E_{22}&0&0\\
			n_3&0&0&0&0\\
			\widetilde{n}_4&0&0&0&0\\
		\end{bNiceArray},
\end{equation*}
where
$$
\rank\begin{bmatrix}
\widetilde{E}_{11}&E_{12}\\
E_{21}&E_{22}
\end{bmatrix}=n_1+n_2.
$$
Since $Q^TE=E^TQ\geq0$, we have
\begin{equation*}
\begin{bmatrix}
(\widetilde{E}_{11})^TQ_{11}+(E_{21})^TQ_{21}&(E_{21})^TQ_{22}\\
(E_{12})^TQ_{11}+(E_{22})^TQ_{21}&(E_{22})^TQ_{22}
\end{bmatrix}\geq0.
\end{equation*}
Because
$$
\rank\begin{bmatrix}
\widetilde{E}_{11}&E_{12}\\
E_{21}&E_{22}
\end{bmatrix}=n_1+n_2
~\text{and}~
\rank\begin{bmatrix}
Q_{11}&0\\
Q_{21}&Q_{22}
\end{bmatrix}=n_1+n_2,
$$
it follows that
\begin{equation*}
\begin{aligned}
&\begin{bmatrix}
(\widetilde{E}_{11})^TQ_{11}+(E_{21})^TQ_{21}&(E_{21})^TQ_{22}\\
(E_{12})^TQ_{11}+(E_{22})^TQ_{21}&(E_{22})^TQ_{22}
\end{bmatrix}>0\\
&\Longrightarrow~(E_{22})^TQ_{22}>0\\
&\Longrightarrow~\rank(E_{22})=n_2.
\end{aligned}
\end{equation*}

\textit{Step 6.}
Let
\begin{equation*}
\widetilde{U}A\widetilde{V}=
\begin{bNiceArray}{cccc}[first-row,first-col,nullify-dots]
&n_1&n_2& n_3& \widetilde{n}_4\\
n_1&A_{11}&A_{12}&A_{13}&\widetilde{A}_{14}\\
n_2&A_{21}&A_{22}&A_{23}&\widetilde{A}_{24}\\
n_3&A_{31}&A_{32}&A_{33}&\widetilde{A}_{34}\\
\widetilde{n}_4&\widetilde{A}_{41}&\widetilde{A}_{42}&\widetilde{A}_{43}&\widetilde{A}_{44}\\
\end{bNiceArray}.
\end{equation*}
Compute orthogonal matrices $U_4$ and $V_4$ such that
\begin{equation*}
U_4\widetilde{A}_{44}V_4=
\begin{bNiceArray}{cc}[first-row,first-col,nullify-dots]
&n_4&n_5\\
n_4&A_{44}&0\\
n_5&0&0\\
\end{bNiceArray},
\end{equation*}
where $\rank(A_{44})=n_4$ and $\widetilde{n}_4=n_4+n_5$.
Finally, define
$$
U=\begin{bmatrix}
I&0&0&0\\
0&I&0&0\\
0&0&I&0\\
0&0&0&U_4\\
\end{bmatrix}\widetilde{U},
~
W=\widetilde{W},
~
V=\widetilde{V}\begin{bmatrix}
I&0&0&0\\
0&I&0&0\\
0&0&I&0\\
0&0&0&V_4\\
\end{bmatrix}.
$$
Then the matrices $UEV$, $UBW$, $W^TCV$, and $UAV$ have exactly the required form in Lemma \ref{le202507221916}.~\hfill$~~\square$
\end{document}